\documentclass[a4paper]{amsart}
\usepackage{amssymb,verbatim,dsfont}
\usepackage{xypic}
\xyoption{all}
\usepackage{verbatim}

\newtheorem{theorem}{Theorem}[section]
\newtheorem{lemma}[theorem]{Lemma}
\newtheorem{proposition}[theorem]{Proposition}
\newtheorem{corollary}[theorem]{Corollary}

\theoremstyle{definition}

\theoremstyle{remark}

\newcommand{\ot}{\otimes}

\newcommand{\ov}{\overline}

\begin{document}

\title{A characterization of quiver algebras based on double derivations}


\author{Jorge A. Guccione}
\address{Departamento de Matem\'atica\\ Facultad de Ciencias Exactas y Naturales, Pabell\'on 1
- Ciudad Universitaria\\ (1428) Buenos Aires, Argentina.} \curraddr{} \email{vander@dm.uba.ar}
\thanks{UBACYT X095 and PICT-2006-00836}

\author{Juan J. Guccione}
\address{Departamento de Matem\'atica\\ Facultad de Ciencias Exactas y Naturales\\
Pabell\'on 1 - Ciudad Universitaria\\ (1428) Buenos Aires, Argentina.} \curraddr{}
\email{jjgucci@dm.uba.ar}
\thanks{UBACYT X095 and PICT-2006-00836}

\begin{abstract} Let $k$ a characteristic zero field. We give a characterization for the
finite quiver $k$-algebras, based on double derivations. More precisely, we prove that if an
associative and unitary $k$-algebra have a family of double derivations satisfying suitable
conditions, then it is (canonically isomorphic to) a quiver algebra. This is the
non-commutative version of a result of D. Wright.\end{abstract}

\date{}


\dedicatory{}


\maketitle

\section*{Introduction}
Let $k$ be a characteristic zero field and $A$ a $k$-algebra. We recall from \cite{B},
\cite{C-E-G} and \cite{G-S} that a double derivation of $A$ is a derivation $D\colon A\to A\ot
A$, where $A$ and $A\ot A$ are considered as $A$-bimodules in the standard way. Suppose by a
moment that $A$ is a commutative ring. In~\cite{W} (see also \cite{C}) it was proved that if
there exist elements $x_1,\dots,x_n\in A$ and commutating derivations $D_1,\dots,D_n\colon A\to
A$ such that

\begin{itemize}

\item[] $D_i(x_j)= \delta_{ij}$,  where $\delta$ is the symbol of Kronecker, for $1\le i,j\le n$,

\smallskip

\end{itemize} then $A$ is the polynomial ring $R[x_1,\dots,x_n]$, where $R$ is the ring of
constants of $A$ (for the definition see the next section).

In this note we give a similar characterization for the finite quiver algebras, but based on
double derivations. Namely, we prove that if an algebra $A$ have a family of double derivations
satisfying suitable conditions, then $A$ is (canonically isomorphic to) a quiver algebra. We
obtain this result as a corollary of another one, which is the exact noncommutative version of
the above mentioned result of \cite{W}. Our method of proof consists in adapt to the
noncommutative setting the one given in that paper.

\section{Preliminaries}
In this note $k$ denotes a characteristic $0$ field,  an algebra means an associative and
unitary $k$-algebra and the unadorned tensor product is the tensor product over $k$. Given a
family $D_1\dots, D_n\colon A\to A\ot A$ of double derivations (see the introduction), we say
that an element $a\in A$ is a {\em constant} of $A$ if $D_i(a)=0$ for all $i$. It is immediate
that the constants  form a subring of $A$.

\smallskip

Let $A$ be an algebra and let $s,t\in A$. It is easy to check that the tensor algebra
$T_A(As\ot tA)$ is isomorphic to the algebra with underlying vector space
$$
A\oplus \bigoplus_{j\ge 0}  As\ot (tAs)^{\ot j}\ot tA
$$
and multiplication
$$
(a_1\ot \cdots \ot a_m)(a_{m+1}\ot \cdots \ot a_n) = a_1\ot \cdots \ot a_ma_{m+1}\ot \cdots \ot a_n.
$$
An isomorphism is the homogeneous map defined as the identity map in degree $0$ and by
$$
(a_1s\ot ta'_1)\ot_A\cdots\ot_A (a_ns\ot ta'_n)\mapsto a_1s\ot\alpha_1\ot\cdots\ot\alpha_{n-1}
\ot ta'_n,
$$
where $\alpha_j = ta'_ja_{j+1}s$, in degree greater than $0$. From now on we will use freely
this identification.

\smallskip

Each double derivation $D$ of $A$  extends to a derivation
$$
D\colon T_A(A\ot A)\to T_A(A\ot A).
$$
via
$$
D(a_0\ot\cdots\ot a_n) = \sum_{i=0}^n a_0\ot\cdots\ot D(a_i)\ot\cdots\ot a_n.
$$
We will say that $D$ is {\em locally nilpotent} if for each $a\in A$ there exists
$n\in\mathds{N}$ (depending on $a$) such that $D^n(a) = 0$.

\begin{proposition} \label{prop 1} If $A$ is an algebra and $D\colon A\to A\ot A$ is a locally
nilpotent double derivation, then the map
$$
\rho\colon A\to T_A(A\ot A),
$$
defined by
$$
\rho(a) = a + D(a) + \frac{D^2(a)}{2!} + \frac{D^3(a)}{3!} + \frac{D^4(a)}{4!} + \cdots
$$
is an injective morphism of algebras.
\end{proposition}

\begin{proof} Straightforward.
\end{proof}

\section{Main result}
Let $A$ be an algebra, $X=\{x_1,\dots, x_n\}$ a set of elements of $A$ and $B=A/\langle
X\rangle$ the quotient algebra of $A$ by the two sided ideal generated by $X$. Given $a\in A$,
we let $\ov{a}$ denote the class of $a$ in $B$.

\begin{theorem} \label{teorema principal} If there exist maps
$$
\xymatrix@R-6.1ex{X \ar[r]^-{s} & A\\
x_i\ar@{|->}[r] & s_i}
\qquad\grow{\xymatrix{{}\save[]+<0pc,-0.3pc>*\txt{and} \restore}}\qquad
\xymatrix@R-6.1ex{X \ar[r]^-{t} & A\\
x_i\ar@{|->}[r] & t_i}
$$
such that $s_ix_it_i = x_i$ for all $i$ and $\{s_1,\dots,s_n,t_1,\dots,t_n\}$ is a set of
idempotent elements of $A$, and there exist locally nilpotent double derivations $D_1,\dots,
D_n$ of $A$ satisfying

\begin{enumerate}

\smallskip

\item $D_i(s_j) = D_i(t_j)=0$ for all $i,j$,

\smallskip

\item $D_i(x_j)=\delta_{ij} s_j\ot t_j$ for all $i,j$,

\smallskip

\item The diagram
$$
\xymatrix@C+8pt{A \rto^-{D_i} \dto^-{D_j} & A\ot A\dto^{D_j\ot A}\\
A\ot A \rto^-{A\ot D_i} & A\ot A\ot A
}
$$
commutes for all $i\neq j$,
\end{enumerate}
then there is an isomorphism $\overline{\rho}\colon A\to T_B(M)$, where
$$
M=\bigoplus_{1\le i \le n} B\ov{s_i}\ot\ov{t_i}B,
$$
such that $\overline{\rho}(x_i) = \ov{s_i}\ot\ov{t_i}$ and $\overline{\rho}^{-1} B $ is the
ring of constants of $A$.
\end{theorem}

\begin{lemma}\label{lema} Let $A$ be an algebra and $D$ a locally nilpotent double
derivation of $A$. Assume that for some $x\in A$ one has $D(x)=s\ot t$, where $s,t\in A$ are
idempotent elements such that $sxt=x$. Assume also that $D(s)=D(t)=0$. Then, the map
$$
\ov{\rho}:=\xymatrix{A\rto^-{\rho} & T_A(A\ot A)\rto^-{\pi} & T_B(B\ov{s}\ot\ov{t}B)},
$$
where $B=A/\langle x \rangle$ and $\pi$ is the canonical surjection, is an isomorphism of
algebras. Moreover, this map identifies $B$ with the ring of constants of $A$.
\end{lemma}

\begin{proof}  Suposse that $\ov{\rho}(a)=0$ for some $a\neq 0$. We assert that for all $n\in
\mathds{N}$ there exists $\sum a^{(1)}\ot\cdots\ot a^{(n)}\in A^{\ot n}$ such that
$$
a=\sum a^{(1)}xa^{(2)}x\cdots xa^{(n-1)}xa^{(n)}.
$$
For $n=1$ this is immediate and for $n=2$ it follows easily from the fact that the class of $a$
in $B$ is zero. Assume that it is true for $n$. Then
$$
\ov{\rho}(a)=0 \Rightarrow \sum \pi D^{n-1}\bigl(a^{(1)}xa^{(2)}\cdots a^{(n-1)}xa^{(n)}
\bigr)=\pi D^{n-1}(a)=0.
$$
Since $D$ is a derivation, $\pi(x)=0$ and $D(x)= x + s\ot t$, this implies that
$$
\sum \pi \bigl(a^{(1)}s\ot t a^{(2)}s\ot \cdots \ot ta^{(n-1)}s\ot t a^{(n)}\bigr) =0.
$$
Hence,
$$
\sum a^{(1)}s\ot t a^{(2)}s\ot \cdots \ot ta^{(n-1)}s\ot t a^{(n)}\in \sum_{i=0}^n A^{\ot i}\ot
AxA\ot A^{\ot n-i},
$$
and the assertion holds for $n+1$. Using now that
\begin{equation}
\sum\rho(a^{(1)})\rho(x)\rho(a^{(2)})\rho(x)\cdots \rho(x)\rho(a^{(n-1)}) \rho(x) \rho(a^{(n)})
= \rho(a)\neq 0,\label{eq1}
\end{equation}
$\rho(x)= x + s\ot t$ and $x=sxt$, it is easy to see that
\begin{equation}
\sum \rho(a^{(1)})s\ot t\rho(a^{(2)})s\ot \cdots \ot t\rho(a^{(n-1)})s\ot t\rho(a^{(n)})\neq
0.\label{eq2}
\end{equation}
In fact, by~\eqref{eq1}, $\rho(a)$ is the sum of the homogeneous terms, obtained by replacing
in the left side of~\eqref{eq2} some of the $(s\ot t)$'s by $sxt = x$. Since, each one of these
terms is the image of the left side of~\eqref{eq2} by an appropriate linear map, if $\rho(a)\ne
0$, then~\eqref{eq2} holds. From~\eqref{eq2} it follows  that the degree of $\rho(a)$ is
greater or equal than $n-1$. Since $n$ is arbitrary this is impossible, and so $\ov{\rho}$ is
an injective map.

Since $\ov{\rho}(x) = \ov{s}\ot \ov{t}$, in order to prove that $\ov{\rho}$ is surjective will
be sufficient to check that its image contains $B$. For $a\in A$, let $n_a$ be the minimal
$n\in \mathds{N}$ such that $D^j(a)\in \ker \pi$ for each $j\ge n$. We are going to prove that
$\ov{a}\in\ov{\rho}A$, by induction on $n_a$. If $n_a=1$, then $\overline{\rho}(a)=\ov{a}$.
Suppose that $n_a=n+1$ and $\ov{c}\in\ov{\rho} A$, for all $c\in A$ with $n_c\le n$. Write
$$
D^n(a)=\sum a_{(1)}\ot \cdots \ot a_{(n+1)}.
$$
(Note that this is the Sweedler notation for the $n$-fold comultiplication of $a$ in a
bialgebra \cite{S}, but here $D$ is a derivation instead of an algebra map). Let
$$
L= \{c\in T_A (A\ot A): \pi D^i(c)=0 \text{ for all $i\ge 0$}\}.
$$
It is easy to see that
\begin{align*}
& D^n(a) - \sum a_{(1)}s\ot t a_{(2)}s\ot \cdots \ot t a_{(n)}s \ot t a_{(n+1)}
\intertext{and}
& \sum D^n\bigl(a_{(1)}xa_{(2)}\cdots a_{(n)}xa_{(n+1)}\bigr) - \sum a_{(1)}s\ot t a_{(2)}s\ot
\cdots \ot t a_{(n)}s \ot t a_{(n+1)}
\end{align*}
belong to $L$. Hence,
$$
\pi D^j\bigl(a - \sum a_{(1)}xa_{(2)}\cdots a_{(n)}xa_{(n+1)}\bigr) = 0\quad\text{for all $j\ge
n$,}
$$
and so, by the inductive hypothesis,
$$
\ov{a - \sum a_{(1)}xa_{(2)}\cdots a_{(n)}xa_{(n+1)}} \in \overline{\rho} A.
$$
Since $\ov{a} = \ov{a - \sum a_{(1)}xa_{(2)}\cdots a_{(n)}xa_{(n+1)}}$ this complete the proof.
\end{proof}

\subsection*{Proof of Theorem~\ref{teorema principal}} For $n=1$ this is Lemma~\ref{lema}.
Suppose that $n>1$ and the result is valid for $n-1$. Again by Lemma~\ref{lema}, we can assume
that
$$
A=T_{B_n} (B_ns_n\ot t_nB_n),
$$
where $B_n=\{a\in A: D_n(a)=0\}$ and $D_n(s_n\ot t_n)=s_n\ot t_n\in B_n\ot B_n$. From~(1) and
~(2) it follows that $s_i,t_i,x_i\in B_n$ for all $i<n$, and from (3), that $D_i B_n\subseteq
B_n$ for all $i<n$. Thanks the inductive hypothesis we can assume that
$$
B_n=T_B(M),
$$
where $B$ is the ring of constants of $A$ for $D_1, \dots, D_n$ and $M=\bigoplus_{j=1}^{n-1}
Bs_j\ot t_jB$. Hence
$$
A=T_{B_n}(B_ns_n\ot t_nB_n)=T_B(M \oplus Bs_n\ot t_nB),
$$
as desired.\qed

\begin{corollary} Under the hypothesis of Theorem~\ref{teorema principal}, if
$\{s_1,\dots,s_n,t_1,\dots,t_n\}$ is a complete set of orthogonal idempotents of the ring of
constants $\overline{\rho}^{-1}B$ of $A$ and
$$
\overline{\rho}^{-1}B=\bigoplus_{e\in \{s_1,\dots,s_n,t_1,\dots,t_n\}} k\, e,
$$
then $A$ is the quiver algebra $kQ$, where $Q$ has vertices $Q_0= \{s_1,\dots,
s_n,t_1,\dots,t_n\}$ and one arrow $s_i \to t_i$ for each $i$.
\end{corollary}

\begin{proof} It is immediate.
\end{proof}

The converse of Theorem~\ref{teorema principal} is also true. That is,  if $B$ is an algebra
and $s_1,\dots,s_n,t_1,\dots,t_n$ are (non-necessarily different) idempotents of $B$, then the
algebra $A = T_B(M)$, where $M = \bigoplus_{i=1}^n Bs_i\ot t_iB$, has double derivations
$D_1,\dots,D_n$ satisfying $(1)--(3)$, and moreover $B$ is the ring of constants of $A$. When
$A$ is the non-commutative polynomial ring $k\{x_1,\dots,x_n\}$ (case $s_1=\cdots=s_n=t_1=
\cdots=t_n=1$) these are the partial double derivations considered in \cite{C-E-G}[Subsection
2.4]. When is a general quiver algebra they are those considered in \cite{B}.

\end{document}